\documentclass[10pt,twoside]{article}

\usepackage[T2A]{fontenc}
\usepackage[cp1251]{inputenc}
\usepackage[english,russian]{babel}
\usepackage[tbtags]{amsmath}
\usepackage{amsfonts,amssymb}

\usepackage{mathrsfs}

\usepackage{graphicx}

\newtheorem{lemma}{Лемма}
\newtheorem{corollary}{Следствие}

\newtheorem{proof}{Доказательство}

\begin{document}


\date{}

\author{Ф.\,А.~Ивлев, А.\,Я.~Канель-Белов, MIPT, BIU, SZU}


\title{Оценка расстояния между двумя телами внутри $n$-мерного шара единичного объёма}

\markboth{Ф.\,A.~Ивлев, А.Я.Белов}{Расстояние между телами внутри единичного шара}

\maketitle

\begin{abstract}
Рассматривается задача об оценке расстояния между двумя телами объёма~$\varepsilon$, расположенными внутри $n$-мерного шара~$U$ единичного объёма, при $n \to \infty$.

Пусть $A$~--- замкнутое множество с гладкой границей объёма $\varepsilon$ ($0 \leq \varepsilon \leq 1/2$), находящееся внутри $n$-мерного шара $U$ единичного объёма, реализующее среди всех множеств объёма $\varepsilon$ множество с минимальной возможной площадью свободной поверхности, лежащих в одном полупространстве относительно некоторой гиперплоскости, прохощящей через центр шара. Тогда $A$ имеет такую же площадь свободной поверхности, как и множество, предствляющее из себя пересечение некоего шара, перпендикулярного $U$, и самого шара $U$.

Библиография: 2 названия.

\end{abstract}

{\bf keywords}
минимальная поверхность, многомерная геометрия, минимальные поверхности, предельные теоремы.

\footnotetext{Работа была проведена с помощью Российского Научного Фонда Грант N 17-11-01377}

\section{Введение и постановка задачи}
\label{sec1}

Пусть $B_n$~--- тело единичного объёма в $n$-мерном пространстве.
При $n \to \infty$ диаметр~$B_n$ тоже стремится к бесконечности,  так что внутри тела~$B$ можно найти далёкие точки.
Это обстоятельство связано с трудностями при переносе конечномерных результатов на бесконечномерные, в частности, в функциональном анализе.
Тем не менее в ряде случае есть основания предполагать, что если взять два множества объёма~$\varepsilon$, то расстояние между ними окажется ограниченной функций от $\varepsilon$, вне зависимости от~$n$. Это обстоятельство может оказаться полезным, в том числе и для переноса результатов на бесконечномерный случай. Интересно, что результаты затрагивают теорию минимальных поверхностей [1], [2].

%


Следующие две гипотезы были предложены Н.~А.~Бобылевым и А.~Я.~Канелем:

\textbf{Гипотеза 1. (случай куба)}
Пусть $\varepsilon$~--- данное число в интервале $(0,1)$, $U_n$~--- $n$-мерный куб единичного объёма.
Внутри $U_n$ выбраны два множества $A$ и $B$, каждое объёма $\varepsilon$.
Тогда расстояние между $A$ и $B$ не больше, чем некоторая константа $D = D(\varepsilon)$, и не зависит от размерности пространства~$n$.

\textbf{Гипотеза 2. (случай шара)}
Пусть $\varepsilon$~--- данное число в интервале $(0,1)$, $U_n$~--- $n$-мерный шар единичного объёма (далее иногда просто $U$).
Внутри $U_n$ выбраны два множества $A$ и $B$, каждое объёма $\varepsilon$.
Тогда расстояние между $A$ и $B$ не больше, чем некоторая константа $D = D(\varepsilon)$, и не зависит от размерности пространства~$n$.

Очевидно, что, если $\varepsilon \geq 1/2$, то обе гипотезы верны.
Поэтому в дальнейшем мы будем предполагать, что  $\varepsilon < 1/2$.

В данной работе доказывается ключевая лемма по нижней оценке площади свободной поверхности множества, находящего внутри $n$-мерного шара единичного объёма, необходимая для доказательства гипотезы 2 в случае произвольных тел.
Продвижение по задаче в случае выпуклых тел $A$ и $B$, находящихся внутри единичного гиперкуба или шара единичного объёма было приведено в предыдущей работе.
%
%
%

\section{Общий ход решения}
\label{subsec2}

{\small
Везде далее будем считать, что рассматриваемые тела гладкие на столько, на сколько это надо. Так как для любого множества в его $\delta$ окрестности можно найти сколько угодно гладкое тело, а замена наших двух тел на их $\delta$ окрестности не сильно влияет на расстояние между ними при $\delta \to 0$, то решив задачу для таких тел, мы получим доказательство и для общего случая.
Определим \textit{площадь свободной поверхности} тела, находящегося внутри шара единичного объёма, как площадь поверхности этого тела без учёта его площади поверхности, являющейся также поверхностью шара $U$, внутри которого решается поставленная задача.

\begin{lemma}
Площадь свободной поверхности тела объёма $\varepsilon$, находящегося внутри $U_n$ не меньше, чем некоторая величина $M(\varepsilon, n)$, равная площади свободной поверхности множества, являющимся пересечением шара $U_n$ с другим шаром, перпендикулярным $U_n$.
\end{lemma}

Доказательство этой леммы будет приведено далее.

Теперь будем расширять наши множества, находящиеся внутри шара внутри этого шара, беря их окрестности. Тогда их объём будет расти пропорционально их площади свободной поверхности (если рассматривать малые приращения), а так как у нас есть на неё оценка снизу, то на уже конечном расширении наших множеств мы получим объёмы равные $1/2$ у каждого из двух множеств, а следовательно, пересечение наших множеств. Тогда расстояние между множествами можно будет оценить как удвоенное необходимое расширение наших множеств до объёма $1/2$, что и завершает доказательство гипотезы.
}

\section{Доказательство ключевой леммы}\label{keyLemma}

\begin{lemma}
	Пусть $A$~--- замкнутое множество с гладкой границей объёма $\varepsilon$ ($0 \leq \varepsilon \leq 1/2$), находящееся внутри $n$-мерного шара $U$ единичного объёма, реализующее среди всех множеств объёма $\varepsilon$ множество с минимальной возможной площадью свободной поверхности, лежащих в одном полупространстве относительно некоторой гиперплоскости, прохощящей через центр шара. Тогда $A$ имеет такую же площадь свободной поверхности, как и множество, предствляющее из себя пересечение некоего шара, перпендикулярного $U$, и самого шара $U$.
\end{lemma}

	Иногда в дальнейшем мы будем использовать выражения ``наше тело'' имея в виду ``множество $A$''. Центр шара $U$ обозначим через $O$. Так же под площадью поверхности тела $A$ почти всегда будет иметься в виду площадь свободной поверхности $A$. Также иногда вместо гиперплоскости мы будем говорить просто плоскость, где понятно, что имеет в виду.

\begin{proof}

Для того, чтобы показать, что наше тело является эквивалентно указанному, мы хотим показать, что наше тело эквивалентно телу вращения. Для того, чтобы это показать мы хотим показать, что оно симметрично относительно некоторого семейства гиперплоскостей, проходящих через центр шара $U$ и перпендикулярных поверхности множеста $A$.

\begin{lemma}\label{div2}
Любое сечение гиперплоскостью, проходящей через центр шара, делящее объём нашего тела пополам, также делит его площадь пополам.
\end{lemma}

\begin{proof}
Если это не так, то можно выбрать ту часть тела, у которой площадь поверхности меньше, и рассмотреть тело, являющееся объединением этой части и части, ей симметричной относительно нашей гиперплоскости. Тогда объём полученного тела будет в точности $\varepsilon$, а площадь поверхности будет не больше, чем удвоенная площадь меньшей из частей, то есть она будет меньше, чем у исходного тела $A$, что противоречит определнию выбора $A$. Значит, такого не может быть.
\end{proof}

\begin{lemma}\label{perp}
Любое сечение из Леммы \ref{div2} в любой точке пересечения с $A$ перпендикулярно поверхности $A$. Или, что то же самое, содержит вектор нормали к поверхности, проведённый в каждой точке пересечения гиперплоскости с нашим телом.
\end{lemma}

\begin{proof}
{\small Приведём здесь набросок доказательства. Если в какой-то точке наша гиперплоскость будет неперпендикулярна поверхности нашего тела, то можно сделать симметрию относительно этой гиперплоскости и получить тело, в котором в окрестности выбранной точки будет ``горбик'': множество точек с бесконечными радиусом кривизны. Но тогда очевидно немного сгладив этот горбик можно при сохранении объёма тела (локально) уменьшить его площадь поверхности (локально), что противоречит определнию $A$.
}
\end{proof}

\begin{lemma}
Пусть $P$~--- точка поверхности. Рассмотрим линейное подпространство, содержащее прямую $OP$ и коразмерности $2$ (в трёхмерии это означает, что оно просто содержит $OP$). Обозначим его через $L$. Будем вращать гиперплоскость, содержащую $L$. а) Тогда существует как минимум одна такая гиперплоскость, которая делит объём и площадь нашего тела пополам. б) Эта гпиерплоскость проходит через нормальный вектор к $A$ в точке $P$ и, следовательно, единственна в случае, когда $L$ не содержит этот вектор нормали. в) Любая гиперплоскость, содержащая вектор нормали и вектор $OP$ делит пополам и площадь и объём множества $A$.
\end{lemma}

\begin{proof}
Существование как минимум одной гиперплоскости вытекает из теоремы о промежуточном вращении, поскольку при фиксации $L$ у нас можно задать функцию ориентированной разности объёмов частей тела $A$, находящихся по разные стороны от вращаемой гиперплоскости и в силу непрерывности этой функции на окружности в коразмерности $2$ она будет иметь 0. А по лемме \ref{div2} оно будет делить пополам не только объём, но и его площадь.
Пункт б) леммы вытекает из леммы \ref{perp}, применённой к пункту а).
Пункт в) следует из того, что если взять в качестве $L$ пространство, получаемое как копространство рассматриваемой гиперплоскости для вектора нормали (гиперплоскость имеет корзамерность 1, а после взятия в ней пространства коразмерности 1, будет как раз пространство коразмерности 2), то в нём по пункту а) такая гиперплоскость должна найтись, но по пункту б) это будет как раз исходная гиперплоскость, потому что есть только одна гиперплоскость, содержащая взятое $L$ и вектор нормали в точке $P$.

\end{proof}

{\small
В дальнейшем будем обозначать данное семейство гиперплоскостей через $\mathcal{A}$. А именно, множество гиперплоскостей, проходящих через $O$ и делящих объём и площадь свободной поверхности множества $A$ пополам.
Рассмотрим процесс симмтерирования: процесс, при котором мы берём и заменяем наше множество на множество являющееся объединением одной из частей нашего тела, на которое поделила его гиперплоскость из семейства $\mathcal{A}$, с часть, полученной из первой путём симметрии относительно этой гиперплоскости.
Заметим, что при процессе симметрирования полученное множество будет обладать всеми теми же самыми свойствами, которые мы доказывали в предыдущих леммах. Действительно, ведь при симметрированнии мы получаем тело того же объёма и той же площади поверхности, что и $A$, то полученное тело так же лежит в классе оптимальных тел для поставленной задачи. Значит, для него верны предыдущие леммы.
}

\begin{lemma}
Семейство $\mathcal{A}$, инвариантно относительно симметрий относительно гиперплоскостей из $\mathcal{A}$.
\end{lemma}

\begin{proof}{\small
Рассмотрим две какие-то гиперплоскости $\alpha$ и $\beta$ из $\mathcal{A}$. Сделаем симметризацию относительно $\alpha$.  Рассмотрим плоскость $\beta$ и какую-нибудь точку $P$ на её пересечении с границей $A$ в той части, которую симметрично отражали в процессе рассматриваемой симметризаци. По доказанному $\beta$ содержит вектор нормали к $A$, восстановленный в точке $P$. С другой стороны, для нового множества любая плоскость, проходящая через $O$ и содержащая этот вектор, должна принадлежать семейству гиперплоскостей, делящих площадь и объём $A'$ пополам (здесь $A'$~--- тело, полученное из $A$ путём рассматриваемой симметризации). Так как выбор $\beta$ относительно $\alpha$ из изначального семейства был произвольным, то любая плоскость из изначального семейства $\mathcal{A}$ будет принадлежать и новому семейству. Аналогично доказывается и обратное утверждение. Значит, рассматриваемые семейства до и после симметризации совпадают.
}
\end{proof}

{\small
Введём на $\mathcal{A}$ топологию, порождённую топологией $n$-мерного проективного просранства на множестве нормалей к этим гиперплоскостям (так как все наши гиперплоскости проходят черзе $O$, то так можно сделать). Заметим, что так как все наши гиперплоскости содержат вектор нормали к $A$ в каждой точке пересечения с $A$ и прим этом через каждый вектор нормали в произвольной точке $A$ проходит целое $n-2$ мерное подсемейство наших гиперплоскостей, то рассматриваемое семейсто $\mathcal{A}$ будет связным в этой топологии.
}

\begin{corollary}
Любое связное семейство гиперплоскостей, имеющих общую точку и инвариантное относительно симметрий относительно любого своего представителя либо соосно (имеет общую ось), либо содержит все гиперплопкости, проходящие через данную точу.
\end{corollary}

Так как множество $A$ лежит в некотором полупространстве относительно одной из гиперплоскостей, проходящих через $O$, то второй вариант невозможен. Откуда следудет, что семейсвто гиперплоскостей $\mathcal{A}$ соосно. При этом заметим, что из непрерывности и связности и того, что в кажой точке тела $A$ к нему можно провести нормальный вектор, мы будем иместь все возможные такие гиперплоскости, которые содержат эту ось.

\begin{lemma}\label{lemma_o_chetvertinkah}
Любые две плоскости из семейства $\mathcal{A}$, перпенидкулярные друг другу, делят множество $A$ на 4 части равного объёма и площади.
\end{lemma}

\begin{proof} {\small
Обозначим две рассматриваемые плоскости $\alpha$ и $\beta$ и пусть они делят множество на 4 части, которые мы обозначим через $K$, $L$, $M$, $N$ таким образом, что $K$ и $L$ лежат а одном полупространстве относительно $\alpha$, а $K$ и $N$ лежат в одном полупространстве относительно $\beta$. Так как плоскости $\alpha$ и $\beta$ взяты из $\mathcal{A}$, то после симметризации относительно $\alpha$ удвоенный объём $K$ должен быть равен удвоенному объёму $L$ (потому что они симметричны относительно $\beta$). А так как сумма их объёмов равна половине всего объёма $A$ (в силу свойсвта плоскости $\alpha$), то каждый из них по объёма будет равен одной четвёртой объёму множества $A$. Аналогично доказывается про две остальных множества и про площадь поверхности.
}
\end{proof}

\begin{lemma}
Любые две плоскости семейства $\mathcal{A}$, находящиеся под углом $\pi / 2^n$ при натуральном $n>1$ высекают на множестве $A$ часть, площадь поверхности и объём которой соответственно равны $1/2^{n+1}$ от площади поверхности и объёма самого множества $A$. Здесь под площадью поверхности этой части имеется в виду площадь поверхности этой части, являющаяся площадью свободной поверхности множеста $A$. Т.\,е. мы не учитываем ``внутреннюю'' площадь поверхности этой части, находящейся по границе гиперплоскостей, по которым мы её высекаем.
\end{lemma}

\begin{proof}
Доказательство этой леммы проходит аналогично доказательству леммы о четвертинках. Для этого достаточно рассмотреть индукцию по $n$ и воспользоваться при переходе леммой о четвертинках и свойствами гиперплоскостей семейства $\mathcal{A}$.
\end{proof}

{\small При любом $n > 1$ мы можем выбрать кусок $A$, находящийся между двумя плоскостями, составляющими угол $\pi/2^n$ мы можем его раскопировать последовательно симметриями относительно плоскостей, составляющими с двумя выбранными углы $k\cdot \pi/2^n$ и получить из маленького кусочка полноценное множество $A'$ того же объёма и той же площади поверхности, поскольку изначальный кусок имел площадь и объёма ровно в $2^{n+1}$ раз меньший, чем $A$. Тогда выбрав достаточно большое $n$ и проделав эту операцию мы получим множество, площадь поверхности и объём которого будет сколько угодно точно близки к площади поверхности и объёму тела, полученного путём вращения некторого сечения, высекаемого некоторой гиперплоскостью семейства $\mathcal{A}$, вокруг общей оси всех гиперплоскостей из $\mathcal{A}$. Значит, в силу предельного перехода можно считать, что наше тело $A$ является телом вращения (можно считать в том смысле, что у него будут те же объём и площадь свободной поверхности).
Тогда задача поиска нужного нам тела сводится к вариационной задачи от одной переменной. А именно, если обозначить через $V_n$ площадь поверхности $n$-мерного шара радиуса $1$, то площадь поверхности нашего множства можно будет посчитать проинтегрировав её вдоль оси вращения:
}
$$
V = \int V_{n-1} \cdot r(x)^{n-1} dx + V_0,
$$

где $V_0$~--- остаток множества $A$, ``вжатый'' в грницу шара $U$.
Утверждается, что решением это задачи как раз и является граница шара, перпендикулярного $U$.
\end{proof}

Отсюда мы получаем необходимую нам экспоненциальную оценку снизу на площадь свободной поверхности тела, объёма $\varepsilon$, находящегося внутри $U$, откуда и следует доказательство исходной гипотезы.

MIPT, 

College of Mathematics and Statistics, Shenzhen University, Shenzhen, 518061, China

%
%
%
%
%



\end{document}